\documentclass[titlepage,twoside,12pt]{article}
\usepackage{amssymb}
\usepackage{amsfonts}
\begin{document}

{\LARGE \bf New Symmetry Groups for \\ \\Generalized Solutions of ODEs} \\ \\

{\bf Elem\'{e}r E ~Rosinger} \\ \\
{\small \it Department of Mathematics \\ and Applied Mathematics} \\
{\small \it University of Pretoria} \\
{\small \it Pretoria} \\
{\small \it 0002 South Africa} \\
{\small \it eerosinger@hotmail.com} \\ \\

{\bf Abstract} \\

It is shown for a simple ODE that it has many symmetry groups beyond its usual Lie group
symmetries, when its generalized solutions are considered within the nowhere dense
differential algebra of generalized functions. \\ \\

{\bf 0. Idea and Motivation} \\

The standard Lie group theory applied to symmetries of solutions of PDEs, Olver [1-3], deals
with {\it classical}, and specifically, ${\cal C}^\infty$-smooth solutions of such equations.
On the other hand, as is well known, and especially in the case of nonlinear PDEs, there is a
major interest in solutions which are no longer classical, thus in particular, fail to be
${\cal C}^\infty$-smooth, and instead are {\it generalized solutions}. \\

In Rosinger [2-13] a characterization and construction was given for the {\it infinitely} many
differential algebras of generalized functions, each such algebra containing the Schwartz
distributions. These algebras prove to be particularly suitable, among others, for finding
generalized solutions to large classes of nonlinear PDEs, Rosinger [1,6-11,13,15],
Oberguggenberger. \\
The construction of such algebras overcomes the celebrated and often misunderstood 1954
Schwartz impossibility of multiplying distributions. And in doing so, it shows that the
multiplication of generalized functions, and in particular, of distributions, can {\it
naturally}, and in fact, {\it inevitably} be done in infinitely many different ways. More
precisely, in each of these infinitely many differential algebras of generalized functions,
the multiplication of ${\cal C}^\infty$-smooth functions is the same with their usual
multiplication. However, the multiplication of functions which are less than ${\cal
C}^\infty$-smooth, and in particular, the multiplication of distributions or generalized
functions can depend on the specific algebra in which the respective multiplication is
performed. \\
As pointed out, Rosinger [9, pp. 1-9], at the root of the inevitability of having infinitely
many such algebras one finds a rather simple {\it conflict} between multiplication, derivation
and discontinuities. \\

It is precisely this natural and inevitable algebra dependent {\it infinite branching} of the
multiplication of less than ${\cal C}^\infty$-smooth functions, and specifically,
distributions and generalized functions which gives the possibility to find {\it new} symmetry
groups for generalized solutions of nonlinear PDEs. And as with the multiplications themselves,
such new symmetry groups may depend on the differential algebras of generalized functions to
which the respective generalized solutions belong. \\

In this paper we shall present {\it new symmetry groups} of generalized solutions for one of
the simplest possible nontrivial ODEs, when these solutions are considered in the so called
{\it nowhere dense} differential algebras of generalized functions, algebras introduced and
used in Rosinger [1-14], see also Mallios \& Rosinger [1]. \\ \\

{\bf 1. Review of Lie Group Actions} \\

Given a linear or nonlinear PDE \\

(1.1) $~~~ T ( x, D )~ U ( x ) ~=~ 0,~~~ x \in \Omega $ \\

where $\Omega$ is a nonvoid open subset in $\mathbb{R}^n$, one of the major interests in Lie
groups - according to Lie's original aim - is in the study of the symmetries of the solutions
$U : \Omega \longrightarrow \mathbb{R}$ of (1.1), which therefore, transform them into other
solutions of (1.1). \\\

For that purpose, one takes \\

(1.2) $~~~ M = \Omega \times \mathbb{R} $ \\

and finds the Lie groups $G$ corresponding to (1.1), as well as their actions \\

(1.3) $~~~ G \times M \ni ( g, ( x, u ) ) ~\longmapsto~
                         g ( x, u ) \in M,~~~~~ g \in G,~ x \in \Omega,~ u \in \mathbb{R} $ \\

which actions, when extended to the solutions $U \in {\cal C}^\infty ( \Omega, \mathbb{R} )$
of (1.1), will transform them into solutions of the same equation. \\

Here of course in the standard theory, Bluman \& Kumei, Olver [1-3], one only deals with {\it
classical}, thus {\it not} generalized solutions. And in fact, all the way one assumes the
${\cal C}^\infty$-smoothness of solutions. \\

The fact however is that, as is well known, large classes of important solutions of a whole
variety of linear, and especially nonlinear PDEs of interest fail to be classical. And then
the question arises :

\begin{quote}

Whether the standard Lie theory of symmetry of classical, that is, ${\cal C}^\infty$-smooth
solutions can be extended to {\it generalized solutions} as well ?

\end{quote}

A first comprehensive {\it affirmative} answer to that question was presented in Rosinger [1].
There, it was shown that for very large classes of linear and nonlinear PDEs in (1.1) and
their generalized solutions, the following property holds :

\begin{itemize}

\item Whenever a Lie group transforms classical solutions of such an equation into other
solutions of the same equation, that Lie group will also transform the generalized solutions
of that equation into other generalized solutions of that equation.

\end{itemize}

This positive result was obtained by using only one among the infinitely many possible
differential algebras of generalized functions. Namely, all generalized solutions were
elements of the {\it nowhere dense} algebras, the algebras first introduced and used in
Rosinger [1-12,14]. \\

Following the above positive result, a second question arises. Namely, there is a well known
significant difference between the nature of classical, and on the other hand, generalized
solutions, be they solutions of linear or nonlinear PDEs. And then one may naturally ask the
question which, so far, has not been dealt with in the literature :

\begin{quote}

Given a linear or nonlinear PDE in (1.1), are there {\it other} groups associated with it,
beyond the classical Lie groups, and which transform certain generalized solutions of that
equation into other solutions of that equation ?

\end{quote}

In this paper - based on the mentioned infinite branching of the multiplication of generalized
functions - we shall give an {\it affirmative} answer to that question. \\
Here we recall that, as seen in Rosinger [1-16], that infinite branching is already manifested
in differential algebras of generalized functions on $\mathbb{R}$, that is, in the case of one
single independent variable. Therefore, in finding new symmetry groups one may deal with ODEs,
instead of PDEs. \\
Further, for convenience, one may also limit oneself to {\it projectable} Lie group actions,
Olver [1-3], and their corresponding generalizations. \\

The full study of new symmetry groups for generalized solutions of PDEs, and the study of
arbitrary, and not only projectable, such groups is to be undertaken elsewhere. \\ \\

{\bf 2. Projectable Lie Group Actions} \\

Let us consider in more detail the Lie group actions (1.1), (1.2), namely \\

(2.1) $~~~ \begin{array}{l}
                    ~~~G \times M ~\longrightarrow~ M  \\
                   (g, (x,u)) ~\longmapsto~ (\widetilde x,\widetilde u) =
                                                   g(x,u) = (g_1(x,u), g_2(x,u))
           \end{array} $ \\

where $x \in\Omega,~ u \in\mathbb{R}$ are the initial independent and dependent variables,
respectively, while $\widetilde{x} \in\Omega,~ \widetilde{u} \in\mathbb{R}$ are, respectively,
the transformed independent and dependent variables. In other words \\

(2.2) $~~~ \begin{array}{l}
                   G \times M \ni (g,(x,u)) ~\longmapsto~ \widetilde{x} =
                                                                g_1(x,u) \in \Omega \\
                   G \times M \ni (g,(x,u)) ~\longmapsto~ \widetilde{u} =
                                                                g_2(x,u) \in \mathbb{R}
           \end{array} $ \\

with $g_1$ and $g_2$ being ${\cal C}^\infty$-smooth. \\

We note that, given $g \in G$, in view of the Lie group axioms, it follows that the mapping \\

(2.3) $~~~ M \ni (x,u) ~\stackrel{g}{\longmapsto}~ g(x,u) \in M $ \\

is a ${\cal C}^\infty$-smooth diffeomorphism. Thus we have the {\it injective group
homomorphism} \\

(2.4) $~~~ G \ni g ~\longmapsto~ f_g \in {\cal D}iff^\infty ( M, M ) $ \\

where $f_g$ is defined by \\

(2.5)~~~ $ M \ni (x,u) ~\longmapsto~ f_g (x,u) = g (x,u) \in M $ \\

Here the noncommutative group structure on ${\cal D}iff^\infty ( M, M )$ is defined by the
usual composition of mappings, and thus the neutral element is $e = id_M$, that is, the
identity mapping of $M$ onto itself. \\

In this way, in terms of the Euclidean domain $M$, the group homomorphism (2.4) is but a group
of smooth coordinate transforms, each of which has a smooth inverse. \\

Now, the Lie group actions (2.1), (2.2) are called {\it projectable}, Olver [1-3], if and only
if $g_1$ has the particular form of {\it not} depending on $u$, namely \\

(2.6) $~~~ \widetilde{x} = g_1(x,u) = g_1(x),~~~ g \in G,~ (x,u) \in M $ \\

The advantage of such projectable Lie group actions (2.6) comes from the fact that, in view of
(2.3), for $g \in G$, we obtain the ${\cal C}^\infty$-smooth diffeomorphism \\

(2.7) $~~~ \Omega \ni x ~\stackrel{g_1}{\longmapsto}~ \widetilde{x} = g_1(x) \in \Omega $ \\

Thus given any function $U : \Omega \longrightarrow \mathbb{R}$, it is easy to define the
group action $\widetilde{U} = g U$ on $U$ for any group element $g \in G$, namely \\

(2.8) $~~~ \widetilde{U} ( x ) ~=~ ( g U ) ( x ) ~=~
                           g_2 ( g^{-1}_1 ( x ), U ( g^{-1}_1 ( x ) ) ),~~~ x \in \Omega $ \\

In this way projectable Lie group actions (2.1), (2.2), (2.6) can easily be extended to group
actions on ${\cal C}^\infty$-smooth functions \\

(2.9) $~~~ G \times {\cal C}^\infty ( \Omega, \mathbb{R} ) \ni ( g, U )
         ~\longmapsto~ \widetilde{U} = g U \in {\cal C}^\infty ( \Omega, \mathbb{R} )$ \\

given by (2.8). \\ \\

{\bf 3. Nowhere Dense Differential Algebras of Generalized \\
\hspace*{0.4cm} Functions} \\

The differential algebras used in this paper, and called the {\it nowhere dense algebras}, are
of the form \\

(3.1) $~~~ A_{nd} ( \Omega ) ~=~ ( {\cal C}^\infty ( \Omega, \mathbb{R} ) )^{\bf N} /
                                                          {\cal I}_{nd}( \Omega ) $ \\

where ${\cal I}_{nd}( \Omega )$ is a specially chosen, so called {\it nowhere dense ideal} in
$({\cal C}^\infty( \Omega, \mathbb{R} ) )^{\bf N}$, see (3.3) below.  \\
Here we denoted by $( {\cal C}^\infty( \Omega, \mathbb{R} ) )^{\bf N}$ the set of all
sequences $s = ( s_0, s_1, s_2, ~.~.~.~ )$ of functions $s_\nu \in {\cal C}^\infty ( \Omega,
\mathbb{R} )$. Clearly, just like ${\cal C}^\infty( \Omega, \mathbb{R} )$, so is $( {\cal
C}^\infty( \Omega, \mathbb{R} ) )^{\bf N}$ an associative, commutative and unital algebra,
when considered with the usual termwise operations on sequences of functions. \\

Let us here briefly comment on the reasons for the choice of the nowhere dense algebras. Such
a comment may indeed be appropriate in view of the mentioned fact that there are infinitely
many differential algebras of generalized functions to choose from. \\

From the start, let us note that, so far only {\it two} main types of such algebras have been
used in a more consistent manner, although several other ones proved to be useful when dealing
for instance with generalized functions of the Dirac delta type, Rosinger [4-8,14]. \\

Namely, the first to be introduced and systematically applied where the nowhere dense algebras,
Rosinger [1-12,14], and recently, their natural extensions given by {\it space-time foam}
algebras, Rosinger [13,15,16], Mallios \& Rosinger [2]. The main feature of these space-time
foam algebras is that they can easily handle not only singularities on closed, nowhere dense
subsets, but also on arbitrarily large subsets, provided that the complementaries of such
singularity sets are dense. For instance, in case $\Omega = \mathbb{R}$, the singularity set
can be that of all irrational numbers, since its complementary, the rational numbers, is dense
in $\mathbb{R}$. In this way, the space-time foam algebras can handle singularity sets which
are dense, and which have a cardinal larger than that of the set of nonsingular points. \\
All these algebras are defined by conditions which only make use of the topology of their
Euclidean domains of definition $\Omega \subseteq \mathbb{R}^n$. In particular, none of these
algebras involve growth conditions in their definitions. \\

Since the mid 1980s, a second type of algebras introduced in Colombeau became popular among a
number of analysts. These algebras make essential use in their definitions of polynomial type
growth conditions. \\

It is well known, Rosinger [1-12,14], that the nowhere dense algebras $A_{nd} ( \Omega )$
contain the Schwartz distributions, and in fact, contain the ${\cal C}^\infty$-smooth
functions as a {\it differential subalgebra}, namely \\

(3.2) $~~~ {\cal C}^\infty ( \Omega, \mathbb{R} )
           ~\subset~ {\cal D}^\prime ( \Omega ) ~\subset~ A_{nd} ( \Omega ) $ \\

The nowhere dense ideals ${\cal I}_{nd}( \Omega )$ are the set of all sequences of ${\cal
C}^\infty$-smooth functions $ w = ( w_0, w_1, w_2, ~.~.~. ) \in ( {\cal C}^\infty ( \Omega
) )^{\bf N}$ which satisfy the condition \\

(3.3) $~~~ \begin{array}{l}
             \exists~~ \Gamma \subset \Omega,~~ \Gamma ~~\mbox{closed, nowhere dense} ~: \\ \\
             \forall~~ x \in \Omega \setminus \Gamma ~: \\ \\
             \exists~~ V \subset \Omega \setminus \Gamma,~~ V ~\mbox{neighbourhood of}~ x,~~
                                                   \mu \in \mathbb{N} ~: \\ \\
             \forall~~ \nu \in \mathbb{N},~ \nu \geq \mu ~: \\ \\
             ~~~~ w_\nu ~=~ 0~~ \mbox{on}~ V
           \end{array} $ \\

We note that with the termwise partial derivation $D^p$, with $p \in \mathbb{N}^n$, of
sequences of functions, we have \\

(3.4) $~~~ D^p {\cal I}_{nd}( \Omega ) ~\subseteq~ {\cal I}_{nd}( \Omega ) $ \\

thus we can define the arbitrary partial derivatives for generalized functions in the nowhere
dense algebras, by \\

(3.5) $~~~ D^p : A_{nd} ( \Omega ) ~\longrightarrow~ A_{nd} ( \Omega ) $ \\

where \\

(3.6) $~~~ A_{nd} ( \Omega ) \ni s + {\cal I}_{nd} ( \Omega ) ~\longmapsto~
                D^p s + {\cal I}_{nd} ( \Omega ) \in A_{nd} ( \Omega ) $ \\

Finally, related to (3.2) we recall that the inclusion of differential algebras \\

(3.7) $~~~ {\cal C}^\infty ( \Omega, \mathbb{R} ) ~\subset~ A_{nd} ( \Omega ) $ \\

takes place according to the mapping \\

(3.8) $~~~  {\cal C}^\infty ( \Omega, \mathbb{R} ) \ni \psi ~\longmapsto~
         ( \psi, \psi, \psi, ~.~.~.~ ) + {\cal I}_{nd} ( \Omega ) \in A_{nd} ( \Omega ) $ \\

It should be noted that the nowhere dense differential algebras of generalized functions
$A_{nd} ( \Omega )$ where the first in the literature to contain the Schwartz distributions,
Rosinger [4-7], and thus to overcome the 1954 Schwartz impossibility. Furthermore, these
algebras proved to be useful in solving large classes of nonlinear PDEs, in abstract
differential geometry, and in the first complete solution to Hilbert's fifth problem,
Rosinger [1, 4-12,14], Mallios \& Rosinger [1]. \\

A main advantage of the nowhere dense differential algebras of generalized functions $A_{nd}
( \Omega )$ comes from the fact that they are {\it not} defined by any sort of growth
conditions, and instead, they only use the topology on the Euclidean domains $\Omega \subseteq
\mathbb{R}^n$. That fact gives them a significant versatility in dealing with large classes of
nonlinear operations and singularities. Further, it clearly distinguishes them from the later
introduced Colombeau algebras which are defined by polynomial type growth conditions. \\
Also, the nowhere dense algebras $A_{nd}( \Omega )$ form a {\it flabby} sheaf, unlike the
Colombeau, or many other algebras of generalized functions, or for that matter, the Schwartz
distributions. And this flabbiness proves to have important applications, Malios \& Rosinger
[1]. \\

In fact, owing to the structure of the {\it nowhere dense ideals} ${\cal I}_ {nd}(\Omega)$,
the nowhere dense algebras simply do {\it not} notice, are {\it not} sensitive to, or shall we
say, {\it jump} over all singularities on closed, nowhere dense subsets of their domain of
definition. And it should be noted that such a property in handling singularities is
nontrivial, since closed nowhere dense subsets can have {\it arbitrary} large positive
Lebesgue measure, Oxtoby. \\
Clearly, the case of such singularities on subsets $\Gamma$ of positive Lebesgue measure
cannot be treated by the Schwartz distributions, and in particular, by Sobolev spaces. Equally,
they cannot be treated by the generalized functions in the Colombeau algebras. \\

One of the immediate consequences of this treatment of any closed and nowhere dense
singularity is the development in Mallios \& Rosinger [1], which allows a far reaching
extension of the de Rham cohomology, with the incorporation of singularities situated on the
mentioned kind of arbitrary closed and nowhere dense subsets in Euclidean spaces. \\

An earlier consequence of this treatment of singularities was the {\it global}
Cauchy-Kovalevskaia theorem, Rosinger [8-10]. According to that theorem, arbitrary nonlinear
systems of analytic PDEs have global generalized solutions on the whole of their domain of
definition, solutions given by elements of the nowhere dense algebras. Furthermore, these
global generalized solutions are usual analytic functions on the whole of their domains of
definition, except for closed nowhere dense subsets, which in addition, and if desired, can be
chosen so as to have zero Lebesgue measure. \\

Such a type of very general nonlinear and global existence result has not been obtained in the
Colombeau algebras of generalized functions, owing among others to the {\it growth conditions}
which appear essentially in the definition of those algebras. \\

In Rosinger [7], a wide ranging and purely algebraic, namely, ring theoretic characterization
was given for the first time for all possible differential algebras of generalized functions
which, as in (3.2), contain the Schwartz distributions. \\
Based on that characterization, the Colombeau algebras by necessity were shown to be a
particular case, Rosinger [8,9], Grosser et.al. [p. 7], MR 92d:46098, Zbl. Math. 717 35001, MR
92d:46097, Bull. AMS vol.20, no.1, Jan 1989, 96-101, and MR 89g:35001. \\ \\

{\bf 4. An Example of ODE with Symmetry Groups Beyond Lie \\
\hspace*{0.4cm} Symmetry Groups} \\

Let us consider the simplest nontrivial ODEs which is of the form \\

(4.1) $~~~ U^{\,\prime} ( x ) ~=~ F ( x ),~~~ x \in \Omega = \mathbb{R} $ \\

where $F \in {\cal C}^\infty ( \Omega, \mathbb{R} )$. As is well known, Bluman \& Kumei, the
Lie group symmetry of (4.1) is given by the one dimensional, or one parameter action \\

(4.2) $~~~ G \times M \ni ( \epsilon, ( x, u ) ) ~\longmapsto~ ( x, u + \epsilon ) \in M $ \\

where $G = ( \mathbb{R}, + )$ is the additive group of real numbers, $\epsilon \in G =
\mathbb{R}$ is the one dimensional group parameter, while $M = \Omega \times \mathbb{R} =
\mathbb{R}^2$. \\

This obviously is a projectable Lie group action, thus according to (2.9), it extends easily
to an action on ${\cal C}^\infty$-smooth functions, given by \\

(4.3) $~~~ \mathbb{R} \times {\cal C}^\infty ( \Omega, \mathbb{R} ) \ni ( \epsilon, U )
    ~\longmapsto~ \widetilde{U} = \epsilon U \in {\cal C}^\infty ( \Omega, \mathbb{R} ) $ \\

where \\

(4.4) $~~~ \widetilde{U} ( x ) ~=~ U ( x - \epsilon ),~~~
                       \epsilon \in \mathbb{R},~ x \in \Omega $ \\

We shall now consider the ODE in (4.1) and its solutions - classical and generalized - within
the nowhere dense algebra $A_{nd}( \Omega )$. In this extended context, it turns out that the
ODE in (4.1) has many generalized solutions. Furthermore, these generalized solution admit
many symmetry groups in addition to (4.2). \\

Let us take any function $\rho \in {\cal C}^\infty ( \Omega, \mathbb{R} )$, such that \\

(4.5) $~~~ \rho ~=~ 1 ~~\mbox{on}~~ ( - \infty, -1 ] \cup [ 1, \infty ),~~~
                                      \rho ~=~ 0 ~~\mbox{on}~~ [ - 1/2, 1/2 ] $ \\

Then for any $a, h \in \mathbb{R}$ we define a corresponding action \\

(4.6) $~~~ J_{a,\,h} : ( {\cal C}^\infty ( \Omega, \mathbb{R} ) )^\mathbb{N}
           ~\longrightarrow ~ ( {\cal C}^\infty ( \Omega, \mathbb{R} ) )^\mathbb{N} $ \\

as follows. Given any sequence $s = ( s_0, s_1, s_2, ~.~.~.~ ) \in ( {\cal C}^\infty ( \Omega,
\mathbb{R} ) )^\mathbb{N}$, then \\

(4.7) $~~~ J_{a,\,h}~ s ~=~ ( J_{a,\,h,\,0}~ s_0,~ J_{a,\,h,\,1}~ s_1,~
                                          J_{a,\,h,\,2}~ s_2, ~.~.~.~ ) $ \\

where \\

(4.8) $~~~ J_{a,\,h,\,\nu}~ s_\nu ( x ) ~=~
                 \begin{array}{|l}
                   ~~\rho ( ( \nu + 1 ) ( x - a ) ) s_\nu ( x ) ~~~\mbox{if}~~ x \leq a \\ \\
                   ~~\rho ( ( \nu + 1 ) ( x - a ) ) ( s_\nu ( x ) + h )
                                                               ~~~\mbox{if}~~ x \geq a
                 \end{array} $ \\

for $\nu \in \mathbb{N},~x \in \Omega $ \\

Thus we have for $\nu \in \mathbb{N},~x \in \Omega $ \\

(4.9) $~~~ J_{a,\,h,\,\nu}~ s_\nu ( x ) ~=~
                 \begin{array}{|l}
                    ~~s_\nu ( x ) ~~~\mbox{if}~~ x \leq a - 1 / ( \nu + 1 ) \\ \\
                    ~~0 ~~~\mbox{if}~~ a - 1 / 2 ( \nu + 1 ) \leq x
                                                \leq a + 1 / 2 ( \nu + 1 ) \\ \\
                    ~~s_\nu ( x ) + h ~~~\mbox{if}~~ x \geq a + 1 / ( \nu + 1 )
                 \end{array} $ \\

It follows that we have \\

(4.10) $~~~ \begin{array}{l}
               \forall~~~ t, s \in ( {\cal C}^\infty ( \Omega,
                                         \mathbb{R} ) )^\mathbb{N} ~: \\ \\
               ~~~~ t - s \in {\cal I}_{nd} ( \Omega ) ~~~\Longrightarrow~~~
                               J_{a,\,h}~ t - J_{a,\,h}~ s \in {\cal I}_{nd} ( \Omega )
             \end{array} $ \\

Indeed, according to (3.3) and the hypothesis, there exists a closed and nowhere dense subset
$\Gamma \subset \Omega$, such that for every $x \in \Omega \setminus \Gamma$, there exists a
neigbourhood $V \subset \Omega \setminus \Gamma$ of $x$ on which $t_\nu = s_\nu$, for $\nu \in
\mathbb{N}$ large enough. \\
But $\Gamma_a = \Gamma \cup \{ a \}$ is still closed and nowhere dense in $\Omega$. And it is
easy to see that, for $\Gamma_a$ and $t - s$, the condition (3.3) holds. \\

It follows that we can define the action \\

(4.11) $~~~ J_{a,\,h} : A_{nd} ( \Omega ) ~\longrightarrow~ A_{nd} ( \Omega ) $ \\

by \\

(4.12) $~~~ A_{nd} ( \Omega ) \ni s + {\cal I}_{nd} ( \Omega ) ~\longmapsto~
                J_{a,\,h}~ s + {\cal I}_{nd} ( \Omega ) \in A_{nd} ( \Omega ) $ \\

The main point of interest is the following {\it commutative group} property of the action in
(4.11), (4.12), namely \\

(4.13) $~~~ J_{a,\,h} \circ J_{a,\,k} ~=~ J_{a,~ h + k},~~~ h,~ k \in \mathbb{R} $ \\

Indeed, in view of (4.9), we have for every sequence $s = ( s_0, s_1, s_2, ~.~.~.~ ) \in
( {\cal C}^\infty ( \Omega, \mathbb{R} ) )^\mathbb{N}$ the relation \\

$ ( J_{a,\,h,\,\nu}~ ( J_{a,\,k,\,\nu}~ s_\nu  ) )( x ) ~=~
             \begin{array}{|l}
                ~~s_\nu ( x ) ~~\mbox{if}~ x \leq a - 1 / ( \nu + 1 ) \\ \\
                ~~s_\nu ( x ) + h + k ~~\mbox{if}~ x \geq a + 1 / ( \nu + 1 )
             \end{array} $ \\ \\

Thus obviously $J_{a,\,h,\,\nu}~ ( J_{a,\,k,\,\nu}~ s ) - J_{a,\,h + k,\,\nu}~ s \in {\cal
I}_{nd} ( \Omega )$. \\

Finally, we can return to the ODE in (4.1). Let $U \in {\cal C}^\infty ( \Omega, \mathbb{R} )$
be any classical solution of it. Then according to (3.6), (3.8), the generalized function \\

(4.14) $~~~ W = ( U, U, U, ~.~.~.~ ) + {\cal I}_{nd} ( \Omega ) \in A_{nd} ( \Omega ) $ \\

is a generalized solution of the ODE in (4.1), when this equation is considered in the nowhere
dense algebra $A_{nd} ( \Omega )$. \\

If we take now any $a, h \in \mathbb{R}$ and apply the corresponding action $J_{a,\,h}$ in
(4.11) to $W$, then clearly \\

(4.15) $~~~ J_{a,\,h} W \in A_{nd} ( \Omega ) ~\setminus~
              {\cal C}^\infty ( \Omega, \mathbb{R} ),~~~ h \in \mathbb{R},~ h \neq 0 $ \\

yet it is easy to see that, with the derivation (3.6) in the nowhere dense algebras, we have \\

(4.16) $~~~ D ( J_{a,\,h} W ) ~=~ F $ \\

that is, $J_{a,\,h} W$ is a {\it generalized solution} of the ODE in (4.1) which is {\it not}
classical when $h \neq 0$. \\

Finally, in view of the group property (4.13), it follows that the generalized solutions (4.15)
of the ODE in (4.1) are transformed in generalized solution of the same equation. \\

It is important to note, however, that the transformations which turn generalized solution of
the ODE in (4.1) into other generalized solution of that equation are far more {\it numerous}
than those presented above. Indeed, let \\

(4.17) $~~~ A ~\subset~ \Omega ~~~\mbox{be any discrete set} $ \\

and let be given any mapping $H$ \\

(4.18) $~~~ A \ni a ~\longmapsto~ h_a \in \mathbb{R} $ \\

Then following the above procedure, one can define an action \\

(4.19) $~~~ J_{A,\,H} : A_{nd} ( \Omega ) ~\longrightarrow~ A_{nd} ( \Omega ) $ \\

which will transform generalized solutions of the ODE in (4.1) into other generalized
solutions of that equation. Furthermore, the composition (4.13) will extend in the case of
actions (4.19) as follows. Let $A$ and $B$ be two discrete subsets of $\Omega$ and let $H : A
\longrightarrow \mathbb{R},~ K : B \longrightarrow \mathbb{R}$ any two mappings. Then \\

(4.20) $~~~ J_{A,\,H} \circ J_{B,\,K} ~=~ J_{C,\,L} $ \\

where \\

(4.21) $~~~ C ~=~ A \cup B $ \\

and $L : C \longrightarrow \mathbb{R}$, such that \\

(4.22) $~~~ L ( c ) ~=~
                \begin{array}{|l}
                     ~~H ( c ) ~~\mbox{if}~~ c \in A \setminus B \\ \\
                     ~~H ( c ) + K ( c ) ~~\mbox{if}~~ c \in A \cap B \\ \\
                     ~~ K ( c ) ~~\mbox{if}~~ c \in B \setminus A
                \end{array} $ \\


\begin{thebibliography}{99}


\bibitem{} Bluman G W, Kumei S : Symmetries and Differential Equations. Applied Mathematical
Sciences Vol. 81, Springer, New York, 1989

\bibitem{} Colombeau J F : New Generalized Functions and Multiplication of Distributions.
Mathematics Studies, Vol. 84, North-Holland, Amsterdam

\bibitem{} Grosser M, Kunzinger M, Oberguggenberger M, Steinbauer R : Geometric Theory of
Generalized Functions. Kluwer, Dordrecht, 2002

\bibitem{} Mallios A, Rosinger E E [1] : Abstract differential geometry, differential
algebras of generalized functions, and de Rham cohomology. Acta Applicandae Mathematicae,
Vol. 55, No. 3, February 1999, pp. 231 - 250

\bibitem{} Mallios A, Rosinger E E [2] : Space-time foam dense singularities and de Rham
cohomology. Acta Applicandae Mathematicae, Vol. 67, No. 1, May 2001, 59 - 89,
arXiv:math.DG/0406540

\bibitem{} Oberguggenberger M B : Multiplication of Distributions and Applications to PDEs.
Pitman Research Notes in Mathematics, Vol. 259, Longman, Harlow, 1992

\bibitem{} Olver P J [1] : Applications of Lie Groups to Differential Equations. Springer,
New York, 1986

\bibitem{} Olver P J [2] : Equivalence, Invariants and Symmetry. Camberidge Univ. Press, 1995

\bibitem{} Olver P J [3] : Nonassociative local Lie groups. J. Lie Theory, Vol. 6, 1996, 23-51

\bibitem{} Oxtoby J C : Measure and Category. Springer, New York, 1971

\bibitem{} Rosinger E E [1] : Parametric Lie Group Actions on Global Generalized Solutions of
Nonlinear PDEs, including a Solution to Hilbert's Fifth Problem, (234 pages). Kluwer,
Dordrecht, 1998

\bibitem{} Rosinger E E [2] : Embedding the D$\,^\prime$ Distributions into Pseudotopological
Algebras. Stud. Cerc. Mat. Acad. Romania, Vol. 18, No. 5, 1966, pp. 687-729

\bibitem{} Rosinger E E [3] : Pseudotopological Spaces, Embedding the D$\,^\prime$
Distributions into Algebras. Stud. Cerc. Mat. Acad Romania, Vol. 20, No. 4, 1968, pp. 553-582

\bibitem{} Rosinger E E [4] : Division of Distributions. Pacific Journal of Mathematics,
Vol. 66, No. 1, 1976, pp. 257-263

\bibitem{} Rosinger E E [5] : Nonsymmetric Dirac Distributions in Scattering Theory.
Proceedings of the Fourth Conference on Ordinary and Partial Differential Equations, Dundee,
Scotland, 1976, Springer Lecture Notes in Mathematics, Vol. 584, pp. 391-399

\bibitem{} Rosinger E E [6] : Distributions and Nonlinear Partial Differential Equations.
Springer Lecture Notes in Mathematics, Vol. 684, 1978 (146 pages)

\bibitem{} Rosinger E E [7] : Nonlinear Partial Differential Equations Sequential and Weak
Solutions. North-Holland Mathematics Studies, Vol. 44, 1980 (317 pages)

\bibitem{} Rosinger E E [8] : Generalized Solutions of Nonlinear Partial Differential
Equations. North Holland Mathematics Studies, Vol. 146, 1987 (409 pages)

\bibitem{} Rosinger E E [9] : Nonlinear Partial Differential Equations, An Algebraic View
of Generalized of Generalized Solutions. North-Holland Mathematics Studies, Vol. 164,
1990 (380 pages)

\bibitem{} Rosinger E E [10] : Global version of the Cauchy-Kovalevskaia theorem for nonlinear
PDEs. Acta Appl. Math., Vol. 21, 1990, 331-343

\bibitem{} Rosinger E E [11] : Characterization for the Solvability of Nonlinear PDEs.
Transactions of the American Mathematical Society, Vol. 330, No. 1, March 1992, pp. 203-225

\bibitem{} Rosinger E E [12] : Arbitrary Global Lie Group Actions on Generalized Solutions of
Nonlinear PDEs and an Answer to Hilbert's Fifth Problem. In (Eds. Grosser M, H\"{o}rmann G,
Kunzinger M, Oberguggenberger M B) Nonlinear Theory of Generalized Functions, 251-265,
Research Notes in Mathematics, Chapman \& Hall / CRC, London, New York, 1998

\bibitem{} Rosinger E E [13] : How to solve smooth nonlinear PDEs in algebras of generalized
functions with dense singularities (invited paper) Applicable Analysis, vol. 78, 2001,
355-378, arXiv:math.AP/0406594

\bibitem{} Rosinger E E [14] : Scattering in highly singular potentials.
arXiv:quant-ph/0405172

\bibitem{} Rosinger E E [15] : Space-time foam differential algebras of generalized functions
and a global Cauchy-Kovalevskaia theorem (revised). Technical Report UPWT 99/8, May 1999
(33 pages)

\bibitem{} Rosinger E E [16] : Differential algebras with dense singularities on manifolds.
Technical Report UPWT 99/9 June 1999 (34 pages)

\bibitem{} Rosinger E E, Walus E Y [1] : Group invariance of generalized solutions obtained
through the algebraic method. Nonlinearity, Vol. 7, 1994, 837-859

\bibitem{} Rosinger E E, Walus E Y [2] : Group invariance of global generalised solutions of
nonlinear PDEs in nowhere dense algebras. Lie Groups and their Applications, Vol. 1, No. 1,
July-August 1994, 216-225

\end{thebibliography}
\end{document}